\documentclass{article}
\usepackage{amsmath,amssymb}

\usepackage{tikz}
\usetikzlibrary{calc,decorations.markings}

\usepackage{graphicx}

\let\e\epsilon
\def\C{{\mathbb C}}
\def\S{{\mathbb S}}

\def\cG{{\mathcal G}}

\def\gl{{\mathfrak{gl}}}
\def\sl{{\mathfrak{sl}}}
\def\so{{\mathfrak{so}}}

\def\corank{{\rm corank}}

\newtheorem{theorem}{Theorem}[section]
\newtheorem{lemma}[theorem]{Lemma}
\newtheorem{remark}[theorem]{Remark}
\newtheorem{corollary}[theorem]{Corollary}
\newtheorem{example}[theorem]{Example}

\title{Polynomial graph invariants induced from the $\gl$-weight system}

\author{N.~Kodaneva \thanks{Higher School of Economic},
S.~Lando
\thanks{Higher School of Economics, Skolkovo Institute of Science and Technology
}}

\begin{document}

\maketitle
\date{}

\begin{center}
{\bf Notation}
\end{center}

\begin{tabular}{rl}

$\gl(N)$ & Lie algebra of all matrices $N\times N$\\
$E_{i,j}$ & matrix unit, $1\le i,j\le N$\\
$w_\gl$ & weight system associated to Lie algebra $\gl$\\
$U\gl(\infty)$ & universal enveloping algebra of the limit Lie algebra $\gl(\infty)$\\
$ZU\gl(N)$ & center of $U\gl(N)$\\
$C_1,C_2,\dots$ & Casimir elements generating the center of~$U\gl(\infty)$\\

$\alpha\in\S_m$ & permutation, which determines the hyperedges of the hypermap\\
$\S_m$ & group of permutations of $m$ elements\\
$G(\alpha)$ & the digraph of a permutation $\alpha$\\
$f(\alpha)$ & number of faces of~$\alpha$\\
$c(\alpha)$ & number of cycles in~$\alpha$\\
$\alpha|_U$ & subpermutation of~$\alpha$\\
$\sigma=(1,2,\dots,m)\in\S_m$ & standard long cycle, the hypervertex of the hypermap\\
$\varphi\in\S_m$ & permutation, which determines the hypefaces of the hypermap\\

${\rm St}_N$ & standard representation of $\gl(N)$\\
$F(N)$ & $w_\gl$ in the standard representation\\
$F_\e(N)$ & an $\e$-deformation of the standard representation\\

\end{tabular}

\begin{tabular}{rl}

$G$ & simple graph\\
$K_n$ & complete graph with $n$ vertices\\
$V(G)$ & the set of vertices of~$G$\\
$A_G$ & adjacency matrix of~$G$\\
$G|_U$ ? & the subgraph of $G$ induced by~$U\subset V(G)$\\
$G^{ab}$ & pivot of $G$ on its edge $ab$\\
$G\setminus a$ & graph obtained from~$G$ by erasing the vertex~$a$\\

$\nu(G)$ & nondegeneracy of~$G$\\
$Q_G(u)$ & skew-characteristic polynomial of~$G$\\
$\overline{Q}_G(u,w)$ & refined skew-characteristic polynomial of $G$\\
$L_G(x)$ & interlace polynomial of~$G$\\

$B$ & a chord or an arc diagram\\
$\pi(B)$ & projection of $B$ onto the subspace of primitive elements\\
$\gamma(B)$ & intersection graph of $B$\\
$I(U)$ & the subset of arcs one or both of whose ends belong to the set~$U$ of ends\\

$D=(E;\Phi)$ & delta-matroid\\
$d_D$ & distance in a delta-matroid~$D$\\
$*$ & partial duality in delta-matroids\\

$\cG_n$ & the vector space spanned by graphs with~$n$ vertices\\
$\cG$ & the Hopf algebra of graphs\\
$P(\cG)\subset\cG$ & subspace of primitive elements\\
$D(\cG)\subset\cG$ & subspace of decomposable elements\\
$\pi:\cG\to P(\cG)$ & projection to the subspace of primitive elements\\
& whose kernel is the subspace of decomposable elements in~$\cG$

\end{tabular}


\section{Introduction}
Weight systems are functions on chord diagrams
satisfying so-called Vassiliev's $4$-term relations.
They are closely related to finite type knot invariants,
see~\cite{V90,K93}.

Certain weight systems can be derived
from graph invariants, see a recent account in~\cite{KL22}.
Another main source of weight systems are Lie algebras,
the construction due to D.~Bar-Natan~\cite{BN95}
and M.~Kontsevich~\cite{K93}.
In recent papers~\cite{ZY22,KL22}, the weight systems
associated to Lie algebras~$\gl(N)$ were unified in a
universal $\gl$-weight system, which takes values in
the ring $\C[N,C_1,C_2,C_3,\dots]$ of polynomials in
infinitely many variables. Note that this weight system
is associated to the HOMFLYPT polynomial, which is
an important and powerful knot invariant.
The  unification has been achieved by extending the
$\gl(N)$-weight systems from chord diagrams, which can be considered as involutions without fixed  points modulo
cyclic shifts, to arbitrary permutations.
The main goal of the extension was to produce an efficient
way to compute explicitly the values of the
$\gl(N)$-weight systems, and M.~Kazarian suggested
a recurrence relation for such a computation, which  works
for permutations rather than just for chord diagrams.

A natural question then arises, namely, which  already
known weight systems can be obtained from the universal
$\gl$ weight system.  In addition to understanding the internal
relationship between weight systems, knowing that a
given weight system can be induced from the $\gl$-weight
system  would immediately lead to extending the former
to arbitrary permutations. It is also interesting whether
the universal $\gl$ weight system is related to integrable hierarchies
of partial differential equations in a way similar to umbral
polynomial graph invariants~\cite{KL15,CKL20}.

To each chord diagram, one can associate a graph, called the intersection
graph of the chord diagram.
Certain weight systems are completely determined by the intersection
graphs; this is true, for example, for the $\sl(2)$-
and for the $\gl(1|1)$-weight systems~\cite{CL07}.
This is true no longer for more complicated Lie algebras, say, for $\sl(3)$.
In general, the relationship between Lie algebra weight systems and
polynomial graph invariants looks rather complicated.

Our  main result in the present paper consists in showing
that the well-known graph and delta-matroid invariant,
the interlace polynomial, can be induced from the
universal $\gl$-weight system. We provide an explicit
substitution making the $\gl$-weight system into the
interlace polynomial for chord diagrams and their intersection graphs.
It has been already known that the interlace polynomial
of graphs and delta-matroids  satisfies the $4$-term
relations for them, whence is a weight system~\cite{K20}.
Our proof
yields an independent argument.

We must underline that the property of a graph invariant
to be induced from the $\gl$-weight system is rather
nontrivial.
For example, the chromatic polynomial cannot be induced from the $\gl$-weight system:
the corresponding system of equations for complete graphs $K_1$,  $K_2$, $K_3$, $K_4$, $K_5$
has no solutions.
One should mention, however, the result in~\cite{L00}, where
the chromatic polynomial is realized through a linear combination
of weight systems associated to Lie algebras of the series $\gl$-
and $\so$-, which may hint that it admits an induction from
the universal $\gl$- and $\so$-weight systems
(see~\cite{KY23} for the construction of the latter) treated together.

The results of the present paper lead
to numerous new questions.
For example, it would be interesting to know, whether
the interlace polynomial can be naturally extended to
hypermaps with arbitrarily  many hypervertices  (a
permutation considered modulo  cyclic shifts can be
treated as a hypermap with a single vertex).
Also, there  are universal weight systems
for permutations constructed
through other series of Lie algebras and Lie
superalgebras~\cite{KY23}; their specializations
also can lead to already known, as well as new,
weight systems and their extensions to permutations.

The authors are grateful to M.~Kazarian for numerous
useful discussions and suggestions.

\subsection{$4$-term relations and weight systems}

In V.~Vassiliev theory~\cite{V90} of finite type
knot invariants, a function on chord diagrams
with~$n$ chords is associated to every knot invariant of order not greater than~$n$. He showed that any function constructed
in this way satisfies certain $4$-term relations.
Then in \cite{K93}, M.~Kontsevich showed that every function on chord diagrams satisfying the four-term relations defines a knot invariant of finite type.

Below, we give necessary definitions and state the properties of the objects
we require. A complete account can be found in~\cite{LZ03,CDBook12}.

\begin{figure}
    \centering
    \includegraphics[width=0.8\linewidth]{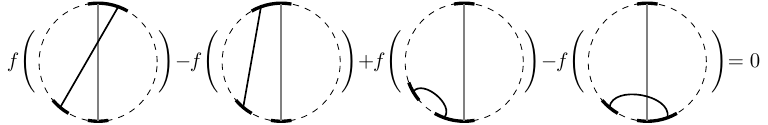}
    \caption{$4$-term relation}
    \label{fig:4t}
\end{figure}

A \emph{chord diagram} of order~$n$ is an oriented circle with $2n$~points on it, which are split into $n$ disjoint~pairs, and each pair is connected by a chord, considered up to orientation
preserving diffeomorphisms of the circle.
A function~$f$ on chord diagrams is called a \emph{weight system} if it satisfies the $4$-term relation shown in Fig.~\ref{fig:4t}. The vector space spanned by chord diagrams
modulo $4$-term relations is endowed with a natural multiplication:
modulo $4$-term relations, the concatenation of two chord
diagrams is well defined.

We can associate the \emph{intersection graph}~$\gamma(B)$ to a given chord diagram~$B$ by associating a vertex to every chord and connecting two vertices by an edge if and only if the
corresponding chords intersect one another. We can also define $4$-term relation for graphs which was introduced in \cite{L00}.

Let $G$ be a graph, and let $a$ and~$b$ be its vertices. Then $4$-term relation is
$$f(G)-f\left(G_{ab}^{\prime}\right)=f\left(\widetilde{G}_{ab}\right)-f\left(\widetilde{G}_{ab}^{\prime}\right),$$
where $G_{ab}^{\prime}$ is obtained from $G$ by deleting the edge between $a$ and~$b$ if it is present in~$G$, or by adding that edge if it is not present in~$G$. And $\widetilde{G}_{ab}$ is obtained from~$G$ by changing the adjacency to~$a$ of all vertices adjacent to~$b$.

\begin{figure}[h]
\begin{center}
\begin{picture}(200,100)(60,0)
\thicklines
\multiput(12,72)(90,0){4}{\circle{40}}
\multiput(55,72)(180,0){2}{{\scriptsize $-$}}
\put(145,72){$=$}
\multiput(5,91)(90,0){4}{\line(-1,-3){10}}
\multiput(-3,84)(90,0){4}{\line(1,0){31}}
\multiput(-8,72)(90,0){4}{\line(3,-1){36}}
\multiput(32,72)(90,0){4}{\line(-5,-2){35}}
\multiput(34,73)(90,0){4}{{\scriptsize b}}
\multiput(12,92)(90,0){4}{\line(1,-5){7.5}}
\put(18,91){\line(-1,-5){7.6}}
\multiput(18,94)(90,0){4}{{\scriptsize a}}
\put(108,91){\line(-2,-3){20.9}}
\put(198,91){\line(1,-2){13}}
\put(288,91){\line(1,-1){12.3}}

\multiput(2,30)(90,0){4}{\circle*{3}}
\multiput(22,30)(90,0){4}{\circle*{3}}
\multiput(2,-10)(90,0){4}{\circle*{3}}
\multiput(22,-10)(90,0){4}{\circle*{3}}
\multiput(-8,10)(90,0){4}{\circle*{3}}
\multiput(32,10)(90,0){4}{\circle*{3}}
\multiput(55,10)(180,0){2}{{\scriptsize $-$}}
\put(145,10){$=$}
\multiput(-6,-16)(90,0){4}{{\scriptsize a}}
\multiput(24,-16)(90,0){4}{{\scriptsize b}}

\multiput(-8,10)(90,0){4}{\line(1,0){39}}
\multiput(-8,10)(90,0){4}{\line(1,2){10}}
\multiput(-8,10)(90,0){2}{\line(1,-2){10}}
\multiput(-8,10)(90,0){4}{\line(3,-2){30}}
\multiput(32,10)(90,0){2}{\line(-3,-2){30}}
\multiput(32,10)(90,0){4}{\line(-1,2){10}}
\multiput(32,10)(90,0){4}{\line(-1,-2){10}}
\multiput(2,30)(90,0){4}{\line(1,0){20}}
\multiput(2,-10)(180,0){2}{\line(1,0){20}}
\multiput(2,-10)(90,0){4}{\line(1,2){20}}

\end{picture}
\end{center}
		\caption{A $4$-term relation for chord diagrams and the corresponding
intersection graphs}
		\label{fourtermrelationfor graphs}
\end{figure}
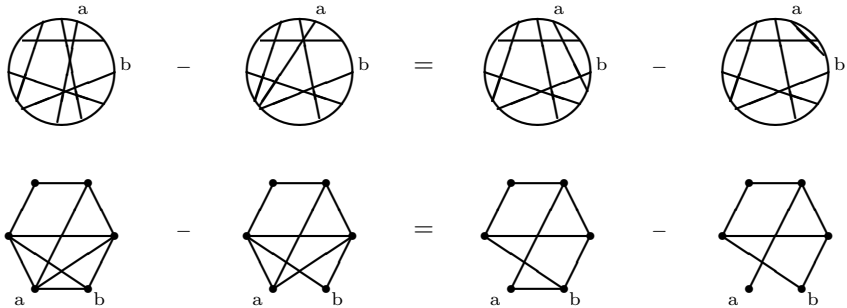

Then every function satisfying the $4$-term relations for graphs defines a weight system and, therefore, a knot invariant of finite type.

\subsection{Graphs, embedded graphs, and their delta-matroids}
Weight systems are closely related to invariants of graphs,
embedded (ribbon, topological) graphs and delta-matroids, see  detailed
account in~\cite{KL22}. Here we recall briefly certain notions
that will be used in the present paper. For an extensive treating of
graph and delta-matroid polynomials in relationship with embedded graphs
see~\cite{EMM13,CMNR19}.

We say that a graph~$G$ is \emph{nondegenerate} if its adjacency matrix~$A_G$ is
nondegenerate mod~$2$.

A chord diagram can be treated as an
orientable embedded graph (a combinatorial map)
 with a single vertex: the vertex is the outer circle, while the edges
(ribbons) are obtained by thickening the chords. It is well known
that the number of faces of the embedded graph associated to a
chord diagram~$B$ coincides with $\corank A_{\gamma(B)}+1$.

A \emph{set system} $D=(E;\Phi)$ is a pair consisting of a finite set~$E$
and a subset $\Phi\subset 2^E$ of the set of its subsets. Below, we consider
only \emph{proper} set systems, i.e. such that $\Phi$ is nonempty.
A set system $D=(E;\Phi)$ is called a \emph{delta-matroid}~\cite{B89,B91,B94} if it
satisfies the following \emph{Symmetric Exchange Axiom}:

for any pair $\phi,\psi\in \Phi$ and any element $x\in\phi\Delta\psi$
there is $y\in\phi\Delta\psi$ such that $\phi\Delta \{x,y\}\in\Phi$.

Here and below $\Delta$ denotes the binary operation of symmetric difference of sets,
$\phi\Delta\psi=(\phi\setminus\psi)\sqcup(\psi\setminus\phi)$.

Elements of~$\Phi$ are called the \emph{admissible sets}
of the delta-matroid~$D$.

Two main sources of delta-matroids are
\begin{itemize}
\item delta-matroids of graphs: for a simple graph~$G$, $D(G)=(V(G);\Phi(G))$,
where a subset $U\subset V(G)$ is admissible iff the induced subgraph $G|_U$
is nondegenerate;
\item delta-matroids of embedded graphs: for a connected embedded graph~$\Gamma$, $D(\Gamma)=(E(\Gamma);\Phi(\Gamma))$,
where a subset $U\subset E(\Gamma)$ is admissible iff the spanning
embedded subgraph $\Gamma|_U$ has a single face.
\end{itemize}

These two constructions are compatible: the delta-matroid $D_{\gamma(B)}$
of the intersection graph $\gamma(B)$ of a chord diagram~$B$ is naturally isomorphic
to the delta-matroid of~$B$ considered as an embedded graph with a single
vertex. The notion of weight system was extended to delta-matroids in~\cite{LZ17}.

\begin{remark}
   The construction of delta-matroids for embedded graphs works both in
   the orientable and non-orientable situations; a similar extension
   for graphs requires considering framed graphs, see~\cite{CMNR18,L06}.
\end{remark}

An important operation on delta-matroids is \emph{partial duality}:
for
a delta-matroid $D=(E;\Phi)$ and a
subset $E'\subset E$ we define $D*E'=(E;\Phi*E')$, where
$\Phi*E'=\{\phi\Delta E'|\phi\in\Phi\}$. If $D=D(\Gamma)$ is the delta-matroid
of an embedded graph~$\Gamma$, then $D*E'$ is the delta-matroid of the
partial dual embedded graph $\Gamma*E'$.

Below, we will require a distance function
for delta-matroids. For a set system $D=(E;\Phi)$, the \emph{distance}
$d_D(U)$ \emph{of a subset}~$U\subset E$ \emph{to}~$D$
is the minimal number of elements in the symmetric difference
of~$U$ with elements of~$\Phi$,
$d_D(U)=\min_{\phi\in\Phi}|U\Delta\phi|$.
The following assertion, which relates the distance
from a subset of edges of an embedded graph
to its delta-matroid with the number of faces
of the embedded graph induced by this subset,
is well-known, but
we were unable to find the explicit statement
in the literature. We decided
to include the proof in the present paper
for completeness.

\begin{lemma}
    Let~$\Gamma$ be an embedded graph, and let
    $D(\Gamma)=(E(\Gamma);\Phi(\Gamma))$
    be its delta-matroid. Then, for a subset
    $E'\subset E(\Gamma)$, we have
    $$
 d_{D(\Gamma)}(E')=f(\Gamma|_{E'})-1,
    $$
where $f(\Gamma_{E'})$ is the number of faces of the
spanning embedded subgraph $\Gamma|_{E'}$.
\end{lemma}

{\bf Proof.} We prove the lemma by induction
on $d_{D(\Gamma)}(E')$. By definition,
this value is~$0$ iff the number of faces
$f(\Gamma||_{E'})$ is~$1$, which provides the
base for the induction.

Now assume that for
any $i$, $0\le i\le k-1$, and any subset $E'\subset E(\Gamma)$ we have
$d_{D(\Gamma)}(E')=i$ iff $f(\Gamma|_{E'})=i+1$.
Pick a subset $E'\subset E(\Gamma)$ such that
$d_{D(\Gamma)}(E')=k$, and let~$\phi\in\Phi(\Gamma)$
be such that $|E'\Delta\phi|=k$. Then for any
element $e\in(\phi\Delta E')$ we have
$d_{D(\Gamma)}(E'\Delta\{e\})=k-1$. Therefore,
by the induction hypothesis, the
number $f(\Gamma|_{'\Delta\{e\}})$ of faces of the embedded spanning subgraph on  the subset
$E'\Delta\{e\}\subset E$ is~$k-1$.
If the  set $E'\Delta\{e\}$ contains the element~$e$,
then by deleting this element we can change the
number of faces of the embedded spanning subgraph
at most by~$1$, whence this number has to become
exactly~$k+1$.  In  the case  the  set $E'\Delta\{e\}$ does not  contain the element~$e$, the similar argument
works for adding the edge~$e$ to the  set $E'\Delta\{e\}$.

Conversely, if  the number of faces of the  embedded
spanning subgraph on the set~$E'$ is~$k+1$, $k>0$, then either this subgraph  is disconnected  or
$E'$  contains an edge~$e$ that bounds two different
faces. In the first case, by adding to~$E'$ an edge
$e\in (E\setminus E')$ connecting two vertices
belonging to two different connected components of
$\Gamma|_{E'}$, we decrease the number of faces by~$1$
making it equal to~$k$. By the induction hypothesis,
there is $\phi\in\Phi(\Gamma)$ such that
$|\phi\Delta(E'\sqcup\{e\})|=k$, so that
$|\phi\Delta E'|\le k+1$, whence equality.
In the second case, by deleting an edge~$e\in E'$
that enters the boundary of two different faces
we also decrease the number of faces by~$1$, and
the argument can be repeated.

$\square$

Permutations considered up to standard cyclic shift can be treated as
orientable hypermaps with a single vertex
A \emph{hypermap} is a triple of permutations $(\sigma,\alpha,\varphi)$ in~$\S_m$ whose product is the identity permutation,
$\sigma\circ\alpha\circ\varphi={\rm id}$, considered up
to simultaneous conjugation by one and the same permutation.
A hypermap is a map if the permutation~$\alpha$ is an involution
without fixed points.
Cycles in the decomposition of $\sigma$ (resp., $\alpha$, $\varphi$) into a product
of independent cycles are \emph{hypervertices}
(resp., \emph{hyperedges}, \emph{hyperfaces}) of the hypermap.
Hypermaps can be considered as either a generalization of maps, or
as bicolored maps~\cite{W75}.

Below, we will always assume that the first permutation
$\sigma$ is the standard cycle, $\sigma=(1,2,\dots,m)$,
and since $\varphi=\alpha^{-1}\circ\sigma^{-1}$ is uniquely
reconstructed from~$\alpha$, we will speak about~$\alpha$
as about a hypermap with a single face.

The number
$f(\alpha)$ of (hyper)faces for a given permutation~$\alpha$
is the number of independent cycles in the permutation~$\varphi$. If~$\alpha$
is an involution without fixed points (an arc diagram), then the number
of its hyperfaces coincides with the  number of ordinary faces.

To each hypermap, an orientable two-dimensional surface
with boundary is naturally associated. The hyperfaces of
the hypermap are in a natural one-to-one correspondence with
the connected components of the  boundary of the surface.
For a single vertex, the surface is constructed by attaching
the hyperedges (which are disks corresponding to the cycles in~$\alpha$)
to the vertex, which is the disk, which corresponds to the unique cycle in~$\sigma$. The rule of gluing is prescribed by~$\alpha$.

\subsection{$\gl$-weight system}
The $\gl$-weight system~$w_\gl$ for permutations is defined in~\cite{ZY22}, see also~\cite{KL22}.
This function on permutations takes values in the algebra $\C[N][C_1,C_2,\dots]$
of polynomials in infinitely many variables $C_1,C_2,\dots$ whose coefficients depend
on the additional variable~$N=C_0$. The variables $C_i$ should be treated
as the standard \emph{Casimir elements}, which generate the center of the universal
enveloping algebra $U\gl(\infty)$ of the limit Lie algebra $\gl(\infty)$.
The definition proceeds as follows.

The value of the $\gl(N)$-weight system on a permutation~$\alpha\in S_m$
is defined as
$$
w_{\gl(N)}(\alpha)=\sum_{i_1,\dots,i_m=1}^NE_{i_1i_{\alpha(1)}}\dots E_{i_mi_{\alpha(m)}},
$$
where $E_{ij}\in\gl(N)$ are the matrix units, $1\le i,j\le N$.

An arc diagram can be interpreted as an involution without fixed points on
the set of the arc ends.
The  $\gl$-weight system for permutations extends
the $\gl(N)$-weight system for
chord diagrams to arbitrary permutations.
This extension provides a way for recursive computation of the $\gl(N)$-weight
systems, as well as unifies the latter weight systems in a universal one, containing all the
specific weight systems, for a given~$N$.

In order to introduce the recurrence, we will need the notion of the digraph of a permutation.
	A permutation can be represented as an oriented graph.
The $m$ vertices of the graph correspond to the permuted elements.
They are placed on the horizontal line, and numbered from left to right
in the increasing order. The arc arrows
show the action of the permutation (so that each vertex is incident with exactly
one incoming and one outgoing arc edge).
	The digraph $G(\alpha)$ of a permutation $\alpha\in \S_m$
consists of these~$m$ vertices and $m$ oriented edges, for example:
	\[
	G(\left(1\ n+1)(2\ n+2)\cdots(n\ 2n)\right)=\begin{tikzpicture}[baseline={([yshift=-.5ex]current bounding box.center)},decoration={markings, mark= at position .55 with {\arrow{stealth}}}]
		\draw[->,thick] (-2,0)--(2,0);
		\draw[blue,postaction={decorate}] (-1.8,0) ..  controls (-1,.5) ..(.2,0);
		\draw[blue,postaction={decorate}] (-1.4,0) ..  controls (-.6,.5) ..(.6,0);
		\draw[blue,postaction={decorate}] (-.2,0) ..  controls (1,.5) ..(1.8,0);
		\draw[blue,postaction={decorate}] (.2,0) ..  controls (-1,-.5) ..(-1.8,0);
		\draw[blue,postaction={decorate}] (.6,0) ..  controls (-.6,-.5) ..(-1.4,0);
		\draw[blue,postaction={decorate}] (1.8,0) ..  controls (1,-.5) ..(-.2,0);
		\fill[black] (-1.8,0) circle (1pt) node[below] {\tiny 1};
		\fill[black] (-1.4,0) circle (1pt) node[below] {\tiny 2};
		\fill[black] (-0.2,0) circle (1pt) node[below] {\tiny n};
		\fill[black] ( .15,0) circle (1pt) node[below] {\tiny n+1};
		\fill[black] ( .6,0) circle (1pt)  node[below]{\tiny n+2};
		\fill[black] ( 1.8,0) circle (1pt) node[below] {\tiny 2n};
		\node[below] at (-.8,0) {$\cdots$};
		\node[below] at (1.2,0) {$\cdots$};
\end{tikzpicture}
	\]

Chord diagrams are permutations of special kind, namely,
involutions without fixed points.
For them, the initial definition coincides with the one above.

Define \emph{Casimir elements} $C_m\in U\gl(N)$, $m=1,2,\dots$:
$$
C_m=w_{\gl(N)}((1,2,\dots,m))=\sum_{i_1,i_2,\dots,i_m=1}^N E_{i_1,i_2}E_{i_2,i_3}\dots E_{i_m,i_1};
$$
associated to the standard cycles $1\mapsto 2\mapsto3\mapsto\dots\mapsto m\mapsto1$.

\begin{theorem}
The center $ZU\gl(N)$ of the universal enveloping algebra $U\gl(N)$ of $\gl(N)$
is identified with the polynomial ring ${\mathbb C}[C_1,\dots,C_N]$.

\end{theorem}

\begin{theorem}[Zhuoke Yang,~\cite{ZY22}]\label{t-gl}
The $w_{\gl(N)}$ invariant of permutations possesses the following properties:
\begin{itemize}
  \item for the empty permutation, the value of $w_{\gl(N)}$ is equal to $1$;
  \item $w_{\gl(N)}$ is multiplicative with respect to concatenation of permutations;
\item {\rm(}\textbf{Recurrence Rule}{\rm)} For the graph of an arbitrary permutation $\alpha$ in~$\S_m$,
and for any two neighboring elements $l,l+1$, of the permuted set $\{1,2,\dots,m\}$, we have
for the values of the $w_{\gl(N)}$ weight system
\begin{equation*}
  \begin{tikzpicture}[baseline={([yshift=-.5ex]current bounding box.center)},decoration={markings, mark= at position .55 with {\arrow{stealth}}}]
    \draw[->,thick] (-1,0) --  (1,0);
    \fill[black] (-.3,0) circle (1pt) node[below] {\tiny l};
    \fill[black] (.3,0) circle (1pt) node[below] {\tiny l+1};
    \draw (-.5,.8) node[left] {a};
    \draw (-.5,-.8) node[left] {b};
    \draw (.5,.8) node[right] {c};
    \draw (.5,-.8) node[right] {d};
    \draw[blue,postaction={decorate}] (-.5,.8) -- (.3,0);
    \draw[blue,postaction={decorate}] (-.3,0) -- (.5,.8);
    \draw[blue,postaction={decorate}] (-.5,-.8) -- (-.3,0);
    \draw[blue,postaction={decorate}] (.3,0) -- (.5,-.8);
  \end{tikzpicture}-
  \begin{tikzpicture}[baseline={([yshift=-.5ex]current bounding box.center)},decoration={markings, mark= at position .55 with {\arrow{stealth}}}]
    \draw[->,thick] (-1,0) --  (1,0);
    \fill[black] (.3,0) circle (1pt) node[below] {\tiny l+1};
    \fill[black] (-.3,0) circle (1pt) node[below] {\tiny l};
    \draw (-.5,.8) node[left] {a};
    \draw (-.5,-.8) node[left] {b};
    \draw (.5,.8) node[right] {c};
    \draw (.5,-.8) node[right] {d};
    \draw[blue,postaction={decorate}] (-.5,.8) -- (-.3,0);
    \draw[blue,postaction={decorate}] (.3,0) -- (.5,.8);
    \draw[blue,postaction={decorate}] (-.5,-.8) -- (.3,0);
    \draw[blue,postaction={decorate}] (-.3,0) -- (.5,-.8);
  \end{tikzpicture}=
  \begin{tikzpicture}[baseline={([yshift=-.5ex]current bounding box.center)},decoration={markings, mark= at position .55 with {\arrow{stealth}}}]
    \draw[->,thick] (-1,0)  -- (1,0);
    \fill[black] (0,0) circle (1pt) node[above] {\tiny l'};
    \draw (-.5,.8) node[left] {a};
    \draw (-.5,-.8) node[left] {b};
    \draw (.5,.8) node[right] {c};
    \draw (.5,-.8) node[right] {d};
    \draw[blue,postaction={decorate}] (-.5,.8) ..controls (0,.4) .. (.5,.8);
    \draw[blue,postaction={decorate}] (-.5,-.8) -- (0,0);
    \draw[blue,postaction={decorate}] (0,0) -- (.5,-.8);
  \end{tikzpicture}-
  \begin{tikzpicture}[baseline={([yshift=-.5ex]current bounding box.center)},decoration={markings, mark= at position .55 with {\arrow{stealth}}}]
    \draw[->,thick] (-1,0) --  (1,0);
    \fill[black] (0,0) circle (1pt) node[below] {\tiny l'};
    \draw (-.5,.8) node[left] {a};
    \draw (-.5,-.8) node[left] {b};
    \draw (.5,.8) node[right] {c};
    \draw (.5,-.8) node[right] {d};
    \draw[blue,postaction={decorate}] (-.5,-.8) ..controls (0,-.4) .. (.5,-.8);
    \draw[blue,postaction={decorate}] (-.5,.8) -- (0,0);
    \draw[blue,postaction={decorate}] (0,0) -- (.5,.8);
  \end{tikzpicture}
\end{equation*}

\end{itemize}
\end{theorem}

For the special case $\alpha(k+1)=k$, the recurrence looks like follows:
\begin{equation*}
	\begin{tikzpicture}[baseline={([yshift=-.5ex]current bounding box.center)},decoration={markings, mark= at position .55 with {\arrow{stealth}}}]
		\draw[->,thick] (-1,0) --  (1,0);
		\fill[black] (-.3,0) circle (1pt) node[below] {\tiny k};
		\fill[black] (.3,0) circle (1pt) node[below] {\tiny k+1};
		\draw (-.5,.8) node[left] {a};
		\draw (.5,.8) node[right] {b};
		\draw[blue,postaction={decorate}] (-.5,.8) -- (.3,0);
		\draw[blue,postaction={decorate}] (-.3,0) -- (.5,.8);
		\draw[blue,postaction={decorate}] (.3,0) ..controls(0,-.3).. (-.3,0);
	\end{tikzpicture}-
	\begin{tikzpicture}[baseline={([yshift=-.5ex]current bounding box.center)},decoration={markings, mark= at position .55 with {\arrow{stealth}}}]
		\draw[->,thick] (-1,0) --  (1,0);
		\fill[black] (.3,0) circle (1pt) node[below] {\tiny k+1};
		\fill[black] (-.3,0) circle (1pt) node[below] {\tiny k};
		\draw (-.5,.8) node[left] {a};
		\draw (.5,.8) node[right] {b};
		\draw[blue,postaction={decorate}] (-.5,.8) -- (-.3,0);
		\draw[blue,postaction={decorate}] (.3,0) -- (.5,.8);
		\draw[blue,postaction={decorate}] (-.3,0) ..controls(0,-.3).. (.3,0);
	\end{tikzpicture}=C_1\times
	\begin{tikzpicture}[baseline={([yshift=-.5ex]current bounding box.center)},decoration={markings, mark= at position .55 with {\arrow{stealth}}}]
		\draw[->,thick] (-1,0)  -- (1,0);
		\draw (-.5,.8) node[left] {a};
		\draw (.5,.8) node[right] {b};
		\draw[blue,postaction={decorate}] (-.5,.8) ..controls (0,.4) .. (.5,.8);
	\end{tikzpicture}-N\times
	\begin{tikzpicture}[baseline={([yshift=-.5ex]current bounding box.center)},decoration={markings, mark= at position .55 with {\arrow{stealth}}}]
		\draw[->,thick] (-1,0) --  (1,0);
		\fill[black] (0,0) circle (1pt) node[above] {\tiny k'};
		\draw (-.5,.8) node[left] {a};
		\draw (.5,.8) node[right] {b};
		\draw[blue,postaction={decorate}] (-.5,.8) -- (0,0);
		\draw[blue,postaction={decorate}] (0,0) -- (.5,.8);
	\end{tikzpicture}
\end{equation*}

\begin{example}
Let us compute the value of $w_\gl$ on the cyclic permutation $(1\ 3\ 2)$ by switching the places of the nodes $2$ and $3$:
\begin{equation*}
  \begin{tikzpicture}[baseline={([yshift=-.5ex]current bounding box.center)},decoration={markings, mark= at position .55 with {\arrow{stealth}}}]
    \draw[->,thick] (-1,0) --  (1,0);
    \fill[black] (-.5,0) circle (1pt) node[below] {\tiny 1};
    \fill[black] (.5,0) circle (1pt) node[below] {\tiny 3};
    \fill[black] (0,0) circle (1pt) node[below] {\tiny 2};
    \draw[blue,postaction={decorate}] (-.5,0) ..controls(0,.3).. (.5,0);
    \draw[blue,postaction={decorate}] (.5,0) ..controls(.25,-.3).. (0,0);
    \draw[blue,postaction={decorate}] (0,0) ..controls(-.25,-.3).. (-.5,0);
    \draw (0,-0.7) node { $(1\ 3\ 2)$};
    \draw (0,0.8) node {};
  \end{tikzpicture}-
  \begin{tikzpicture}[baseline={([yshift=-.5ex]current bounding box.center)},decoration={markings, mark= at position .55 with {\arrow{stealth}}}]
    \draw[->,thick] (-1,0) --  (1,0);
    \fill[black] (-.5,0) circle (1pt) ;
    \fill[black] (.5,0) circle (1pt) ;
    \fill[black] (0,0) circle (1pt) ;
    \draw[blue,postaction={decorate}] (.5,0) ..controls(0,-.3).. (-.5,0);
    \draw[blue,postaction={decorate}] (-.5,0) ..controls(-.25,.3).. (0,0);
    \draw[blue,postaction={decorate}] (0,0) ..controls(.25,.3).. (.5,0);
    \draw (0,-0.7) node { $(1\ 2\ 3)$};
    \draw (0,0.8) node {};
  \end{tikzpicture}=C_1\times
  \begin{tikzpicture}[baseline={([yshift=-.5ex]current bounding box.center)},decoration={markings, mark= at position .55 with {\arrow{stealth}}}]
    \fill[black] (0,0) circle (1pt) ;
    \draw[->,thick] (-1,0)  -- (1,0);
    \draw[blue,postaction={decorate}] (0,0) ..controls (-.5,.7) and (.5,.7)   .. (0,0);
    \draw (0,-0.47) node { $(1)$};
  \end{tikzpicture}-N\times
  \begin{tikzpicture}[baseline={([yshift=-.5ex]current bounding box.center)},decoration={markings, mark= at position .55 with {\arrow{stealth}}}]
    \draw[->,thick] (-1,0) --  (1,0);
    \fill[black] (-.4,0) circle (1pt) node[above] {};
    \fill[black] (.4,0) circle (1pt) node[above] {};
    \draw[blue,postaction={decorate}] (-.4,0) ..controls(0,-.3).. (.4,0);
    \draw[blue,postaction={decorate}] (.4,0) ..controls(0,.3).. (-.4,0);
    \draw (0,0.8) node {};
    \draw (0,-0.7) node { $(1\ 2)$};
  \end{tikzpicture}
\end{equation*}
\begin{eqnarray*}
w_{\gl(N)}((1\ 3\ 2))&=&w_{\gl(N)}((1\ 2\ 3))+C_1\cdot w_{\gl(N)}((1))-N\cdot w_{\gl(N)}((1\ 2))\\
                &=&C_3+C_1^2-NC_2
\end{eqnarray*}
\end{example}

\begin{corollary}
The $\gl(N)$-weight systems, for $N=1,2,\dots$, are combined into
a universal $\gl$-weight system taking values in the ring of
polynomials in infinitely many variables $\C[N;C_1,C_2,\dots]$.

After substituting a given value of~$N$ and an expression of
higher Casimirs $C_{N+1},C_{N+2},\dots$ in terms of the
lower ones $C_1,C_2,\dots,C_N$, this weight system specifies
into the $\gl(N)$-weight system.
\end{corollary}

We say that a function on permutations \emph{is of $\gl$-type}
if it is invariant with respect to the standard shift
and satisfies the same recurrence relation as the $\gl$-weight
system.

\begin{example}
    The number $f(\alpha)$ of faces of a permutation is a function of $\gl$-type. In fact, it is easy to see that
    the number of faces is preserved by replacing the
    first permutation on the left hand side of the recurrence
    with the first one on the right hand side, as well as
    replacing the second permutation on the left with the
    second one one the right.
\end{example}

\subsection{Refined skew-characteristic polynomial}

The \emph{skew-characteristic polynomial} $Q_G(u)$ of a simple graph~$G$ is defined~\cite{DL22}
as
$$
Q_G(u)=\sum_{{U\subset V(G)
\atop U~\text{nondegenerate}}} u^{|V(G)|-|U|}.
$$
If~$G$ is the intersection graph of a chord diagram, then~$Q_G$ coincides with the characteristic
polynomial of the anti-symmetric adjacency matrix associated to the directed intersection graph $\vec{G}$, see~\cite{DL22}.

The \emph{refined skew-characteristic polynomial} $\overline{Q}_G(u,w)$
is
$$
\overline{Q}_G(u,w)=\sum_{{U\subset V(G)}} u^{|V(G)|-|U|}w^{{{\rm corank} A_U}},
$$
so that $Q_G(u)=\overline{Q}_G(u,0)$ for any graph~$G$.

Both the skew-characteristic polynomial and the refined skew-characteristic polynomial
are $4$-invariants, that is, satisfy $4$-term relations for graphs.

The refined skew-characteristic polynomial admits a natural extension to delta-matroids:
for a delta-matroid $D=(E;\Phi)$, we define
$$
\overline{Q}_D(u,w)=\sum_{{U\subset E}} u^{|E|-|U|}w^{{d_D(U)}},
$$

The refined skew-characteristic polynomial thus defined is a $4$-invariant of delta-matroids, see~\cite{LZ17}.

\subsection{Interlace polynomial}

One way to introduce the interlace polynomial is to define it directly for
delta-matroids. For a delta-matroid~$D=(E;\Phi)$, we set
$$
L_D(x)=\sum_{U\subset E}x^{d_D(U)}.
$$

 The following properties of the interlace polynomial follow easily from the definitions.

\begin{theorem}
    The interlace polynomial of delta-matroids is invariant under the
    partial duality.
\end{theorem}

\begin{theorem}
    The interlace polynomial of delta-matroids coincides with the result of substitution
    $u=1,w=x$ into the refined skew-characteristic polynomial,
    $$
    L_D(x)=\overline{Q}(1,x).$$
\end{theorem}

Originally, the interlace polynomial of graphs was defined recursively~\cite{ABS041}.
In order to introduce it, we firstly define the pivot operation for graphs.

Let $G$ be a graph, and let $a$ and $b$ be two vertices connected by an edge. All vertices except~$a$ and~$b$ are divided into four classes:
\begin{enumerate}
    \item vertices adjacent to both $a$ and $b$;
    \item vertices adjacent to $a$ but not to $b$;
    \item vertices adjacent to $b$ but not to $a$;
    \item vertices adjacent neither to $a$ nor to $b$.
\end{enumerate}
The \emph{pivot}~$G^{ab}$ is the graph obtained from~$G$ by removing all the edges connecting vertices from two different classes
among 1--3 above and adding such edges if they are not present in~$G$.
For chord diagrams and delta-matroids, pivot
acts as partial duality with respect to the subset $\{a,b\}$
od the base set~\cite{B89,B91,C09}.
On the delta-matroid $D_G(V(G);\Phi(G))$ of a graph~$G$, the pivot acts as the
partial duality with respect to the subset $\{a,b\}\subset V(G)$.
In particular, the interlace polynomial of graphs is invariant under pivot.

The interlace polynomial $L_G(x)$ is defined by the initial condition

\emph{if there are no edges in the graph~$G$, then $L_G(x)=(x+1)^n$ where $n$ is the number of vertices in~$G$}

and the following recurrence:

\emph{ if vertices~$a$ and~$b$ of the graph~$G$ are connected by an edge, then
\begin{equation}\label{eIPr}
L_G(x)=L_{G\setminus a}(x)+L_{G^{ab}\setminus b}(x),
\end{equation}
where $G\setminus a$ is the graph obtained from~$G$ by erasing the vertex~$a$ and all the edges connecting $a$ with other vertices}.

\begin{example}
    For the graph \begin{tikzpicture}[baseline={([yshift=-.6ex]current bounding box.center)}]
		\draw[-,thin] (0,0)--(.3,0);
            \draw[-,thin] (0,0)--(0,.3);
            \draw[-,thin] (.3,0)--(.3,.3);
            \draw[-,thin] (0,.3)--(.3,.3);
            \draw[-,thin] (0,0)--(.3,.3);
		\fill[black] ( 0,0) circle (1pt) node{};
		\fill[black] ( .3,.3) circle (1pt)  node{};
		\fill[black] ( .3,0) circle (1pt) node{};
		\fill[black] ( 0,.3) circle (1pt)  node{};
\end{tikzpicture}, the interlace polynomial can be computed as
\[
L\begin{tikzpicture}[baseline={([yshift=-.5ex]current bounding box.center)}]
		\draw[-,thin] (0,0)--(.2,0);
            \draw[-,thin] (0,0)--(0,-.2);
            \draw[-,thin] (.2,0)--(.2,-.2);
            \draw[-,thin] (0,-.2)--(.2,-.2);
            \draw[-,thin] (0,0)--(.2,-.2);
		\fill[black] ( 0,0) circle (1pt) node[above]{\tiny a};
		\fill[black] ( .2,-.2) circle (1pt)  node{};
		\fill[black] ( .2,0) circle (1pt) node[above] {\tiny b};
		\fill[black] ( 0,-.2) circle (1pt)  node{};
\end{tikzpicture}
	(x)=
	L\begin{tikzpicture}[baseline={([yshift=-.5ex]current bounding box.center)}]
            \draw[-,thin] (.2,0)--(.2,-.2);
            \draw[-,thin] (0,-.2)--(.2,-.2);
		\fill[black] ( .2,-.2) circle (1pt)  node{};
		\fill[black] ( .2,0) circle (1pt) node[above] {\tiny b};
		\fill[black] ( 0,-.2) circle (1pt)  node{};
\end{tikzpicture}(x)+
L\begin{tikzpicture}[baseline={([yshift=-.5ex]current bounding box.center)}]
            \draw[-,thin] (0,0)--(0,-.2);
            \draw[-,thin] (0,0)--(.2,-.2);
		\fill[black] ( 0,0) circle (1pt) node[above]{\tiny a};
		\fill[black] ( .2,-.2) circle (1pt)  node{};
		\fill[black] ( 0,-.2) circle (1pt)  node{};
  \end{tikzpicture}(x)=2 L\begin{tikzpicture}[baseline={([yshift=-.5ex]current bounding box.center)}]
            \draw[-,thin] (0,-0.2)--(0.4,-0.2);
		\fill[black] ( 0,-0.2) circle (1pt) node[above]{\tiny a};
		\fill[black] ( .2,-0.2) circle (1pt)  node[above] {\tiny b};
		\fill[black] ( 0.4,-0.2) circle (1pt)  node{};
  \end{tikzpicture}(x)=\]
  \[=2L\begin{tikzpicture}[baseline={([yshift=-.5ex]current bounding box.center)}]
            \draw[-,thin] (0,-0.2)--(0.2,-0.2);
		\fill[black] ( 0,-0.2) circle (1pt) node[above]{\tiny b};
		\fill[black] ( .2,-0.2) circle (1pt)  node{};
  \end{tikzpicture}(x)+ 2L\begin{tikzpicture}[baseline={([yshift=-.5ex]current bounding box.center)}]
		\fill[black] ( 0,-0.2) circle (1pt) node[above]{\tiny a};
		\fill[black] ( .2,-0.2) circle (1pt)  node{};
  \end{tikzpicture}(x)=4L\begin{tikzpicture}[baseline={([yshift=0ex]current bounding box.center)}]
		\fill[black] ( 0,0) circle (1pt) node{};
  \end{tikzpicture}(x)+2(x+1)^2
  =4(x+1)+2(x+1)^2=2x^2+8x+6.
  \]
\end{example}

\begin{remark}
    In some papers the interlace polynomial is normalized differently; namely,
    it is defined as $L_G(x-1)$ in our notation; in other words, the initial
    condition for the discrete graph on~$n$ vertices is~$x^n$ rather than $(x+1)^n$.
\end{remark}

The interlace polynomial admits a two-variable extension,
which is defined for delta-matroids as
$$
\overline{L}_D(x,y)=\sum_{{U\subset E}} x^{|U|}y^{{d_D(U)}}.
$$

The following assertion is obvious:

\begin{theorem} We have
$\overline{Q}_D(u,v)=\overline{L}_D(u^{-1},v)\cdot u^{|E|}$.
\end{theorem}

An element $e\in E$ of the base set of a delta-matroid
$(E;\Phi)$ is called a \emph{loop} if it is not  contained
in any admissible set $\phi\in\Phi$; $e$ is a \emph{coloop}
if it is contained in each admissible set.

\begin{theorem}~\cite{ABS04}
    The two-variable interlace polynomial of delta-matroids satisfies the following recurrence. If $e \in E$ is neither a loop nor a coloop, then
$$
\overline{L}_D(x,y)=\overline{L}_{D \backslash e}(x,y)+x\overline{L}_{D * e \backslash e}(x,y).
$$

If $e$ is a coloop, then
$$
\overline{L}_D(x,y)=(x+y) \overline{L}_{D * e \backslash e}(x,y),
$$
while if $u$ is a loop we have
$$
\overline{L}_D(x,y)=(1+x y) \overline{L}_{D \backslash e}(x,y).
$$
\end{theorem}

These recurrences imply similar recurrences for the refined skew-characteristic polynomial. In particular, if $e \in E$ is neither a loop nor a coloop, then
$$
\overline{Q}_D(u,v)=u\overline{Q}_{D \backslash e}(u,v)+\overline{Q}_{D * e \backslash e}(u,v).
$$

\section{A two-variable family induced from the $\gl$-weight system}

Various specializations of the Casimir elements in the $\gl$-universal weight system
also produce weight systems on chord diagrams. In this section, we define
several such specializations,
which together form a family of weight systems depending on two variables,
and formulate their properties.
We associate certain polynomial graph and delta-matroids invariants with this family.
The proofs of the theorems are given in the next section.

\subsection{Definition of the family}
For each positive integer value of~$N$ denote by  ${\rm St}_N:U\gl(N)\to\gl(N)$
the \emph{standard} representation,
which is identity on $\gl(N)\subset U\gl(N)$, of the Lie algebra $\gl(N)$
and its universal enveloping algebra.
Denote by $T:a\mapsto \frac1N {\rm Tr}\ {\rm St}_N(a)$ the mapping taking an arbitrary element $a\in U\gl(N)$
to $\frac1N$ times the trace of ${\rm St}_N(a)$. Note that division of the trace by~$N$
is required, in particular,
 if we wish that the mapping~$F$ takes the product of chord diagrams
 to the product of polynomials, that is, is multiplicative.



It is easy to see that~$T$ acts on the Casimir elements~$C_m$ in the following way:
$$
T:C_m\to N^{m-1}, \qquad m=1,2,\dots.
$$
Indeed, since
$$
w_{\gl(N)}(C_m)=\sum_{i_1,\dots,i_m=1}^N E_{i_1i_2}E_{i_2i_3}\dots E_{i_mi_1},
$$
all the~$N$ diagonal elements in the $N\times N$-matrix ${\rm St}_N(w_{\gl(N)})(C_m)$
are equal to $N^{m-1}$. In particular, since $w_\gl$ takes a permutation to a polynomial
in $C_1,C_2,\dots$ with coefficients in~$N$, the mapping~$T\circ w_\gl$ takes each permutation to
a polynomial in~$N$. In fact, this polynomial proves to be a monomial:

\begin{theorem}\label{t-sr}
The mapping $F=T\circ w_\gl$ takes an arbitrary permutation $\alpha\in \S_m$ to the monomial
$N^{f(\alpha)-1}$, where $f(\alpha)$ is the number of faces of~$\alpha$.
\end{theorem}

Indeed, according to the definition of the number of faces
of a permutation, all the diagonal elements in the matrix
$T(\alpha)$ are $T^{f(\alpha)-1}$.

\begin{remark}\rm
Note that if~$\alpha$ is an involution without fixed points (an arc diagram),
then this assertion is well known, see~\cite{BN95,LZ03}.
\end{remark}

The \emph{$2$-parameter family $F_\e$ of weight systems} for permutations is the
$\e$-perturbation
of the standard representation substitution, $F_0=F$, obtained as the result of the following
substitution in~$w_\gl$:
$$
C_1:= 1+N\e,\quad C_2:= N+2\e+N\e^2,\quad C_3:=N^2+3N\e+3\e^2+N\e^3,\dots,
$$
$$
C_k:=\frac1N((N+\e)^k-(1-N^2)\e^k),\dots
$$

\subsection{Properties of~$F_\e$}

Let~$\alpha\in \S_m$ be a permutation of~$m$ elements $\{1,2,\dots,m\}$, and let~$U$
be a subset of $\{1,2,\dots,m\}$. Then we define the \emph{subpermutation}
$\alpha|_U$ of~$\alpha$ as the permutation of~$U$ taking each element $u\in U$
to the next element of the cycle in~$\alpha$ containing~$u$, which belongs to~$U$.
We consider the permutation $\alpha|_U$ as a permutation of the set~$\{1,2,\dots,|U|\}$,
the renumbering preserving the relative ordering of the elements permuted. If $|U|=k$, then we say that $\alpha|_U$
is a $k$-\emph{subpermutation} of~$\alpha$.

\begin{theorem}\label{t-fe}
For an arbitrary permutation~$\alpha$ of~$m$ elements $\{1,2,\dots,m\}$,
the value $F_\e(\alpha)$ is
$$
F_\e(\alpha)=\sum_{{U\subset \{1,2,\dots,m\}
}}N^{c(\alpha)-c(\alpha|_U)+f(\alpha|_U)-1}\e^{m-|U|},
$$
where the summation is carried over all the~$2^m$ subsets~$U$ of  $\{1,2,\dots,m\}$.
\end{theorem}

Recall that $c(\alpha)$ denotes the number of cycles in the
decomposition of a permutation~$\alpha$ into the product
of independent cycles.

\begin{example}\label{e-ca}
    For a standard cycle of length~$k$, the substitution
$$
C_k=\frac1N((N+\e)^k-(1-N^2)\e^k)=\sum_{i=1}^{k} {k\choose i-1}N^{i-1}\e^{k-i}+N\e^k
$$
agrees with the assertion of the theorem: each $i$-subpermutation of the standard
$k$-cycle is the standard $i$-cycle, there are ${k\choose i}$ such subpermutations,
and the degree of~$N$ in the contribution of
each of it is $1-1+i-1=i-1$.
\end{example}

\begin{example}
    For the complete graph~$K_2$ on two vertices (that is, for the
    permutation $(13)(24)$), the substitution above yields
$$
F_\e((13)(24))=N^2\e^4+4N\e^3+(2N^2+4)\e^2+4N\e+1.
$$
This value also is in accordance with the assertion of
Theorem~\ref{e-ca}, since for $k\ne 2$ all the $k$-element subsets of $\{1,2,3,4\}$ induce
one and the same, up to a cyclic shift, subpermutation.
Two
out of $2$-element subsets of the set $\{1,2,3,4\}$ form the
two ends of a single arc, and the number of boundary components of the corresponding arc diagram is~$2$; the number of cycles is~$1$, whence the degree of~$N$ in the contribution of
these two subpermutations is $2-1+2-1=2$.
Each of the other four $2$-element subsets
consists of two ends of distinct arcs, and the number of boundary components of
the corresponding hypermap is~$1$, whence the degree of~$N$ in the contribution of
these four subpermutations is $2-2+1-1=0$.
\end{example}

Introduce a new variable $v$, which denotes the perturbation of $F_\e$ on the element $C_2$,
$$
F_\e(K_1)=N+v,\qquad v=2\e+N\e^2.
$$

\begin{theorem}\label{t-rsc}
The value of~$F_\e$ on an arbitrary chord diagram, when expressed in the variable~$v$,
coincides with the refined skew characteristic polynomial $\overline{Q}$ of the intersection
graph of this chord diagram,
$$
F_\e(B)=\overline{Q}_{\gamma(B)}(v,N).
$$
In particular, the value of $F_\e$ on an arbitrary chord diagram is a polynomial in~$v$.
\end{theorem}

Now make the following substitution:
\begin{equation}\label{e1}
N:=z^2-1;\quad \e:=\frac1{1-z}.
\end{equation}
Note that this substitution yields $v=1$.

\begin{theorem}\label{t-is}
For an arbitrary chord diagram~$B$ with at least two chords,
the value $F_\e(B)$ under the above substitution coincides with the
interlace polynomial in the variable~$z^2$
of the intersection graph $\gamma(B)$ of~$B$.
\end{theorem}

\begin{corollary}
The interlace polynomial for orientable maps with one vertex and their
intersection graphs can be induced from the $w_\gl$-weight system.
\end{corollary}

\begin{example}\rm
Consider the chord diagram~$K_3$ consisting of three mutually intersecting
chords; the corresponding permutation is $(14)(25)(36)$
and the intersection graph is~$K_3$, the complete graph on three vertices.

The value of $w_\gl$
on this permutation is
$$
w_\gl((14)(25)(36))=C_2^3+3C_1^2C_2-3NC_2^2+2N^2C_2-2NC_1^2.
$$
Hence,
$$
F_\e((14)(25)(36))=N+6\e+15N\e^2+4(2+3N^2)\e^3+3N(4+N^2)\e^4
+6N^2\e^5+N^3\e^6.
$$
Substitution~(\ref{e1}) takes this polynomial to $4(z^2+1)$,
which is the interlace polynomial $L_{K_3}(z^2)$.

\end{example}
\section{Proof of the main results}

In this section, we give a proof of Theorems~\ref{t-fe}, \ref{t-rsc},
and \ref{t-is}.

\subsection{Proof of Theorem~\ref{t-fe}}

Consider the invariant of permutations defined by the sum in the right-hand side
in the assertion of Theorem~\ref{t-fe}. Denote this invariant by~$\widetilde{F}_\e$. We are going
to prove that $\widetilde{F}_\e=F_\e$.
Clearly, the invariant~$\widetilde{F}_\e$ is multiplicative:
for the concatenation of two permutations, the number of cycles is the sum of the numbers
of cycles in the factors, while the number of boundary components in the product
is one less than the sum of the numbers of boundary components in the factors.

Example~\ref{e-ca} provides the base of the assertion of the theorem. Now, we have only
to prove that $\widetilde{F}_\e$ is a $\gl$-type function.
Pick a permutation~$\alpha$
of the set $\{1,2,\dots,m\}$; all the subsets of this set split into three classes:
those containing both elements $l,l+1$ entering the recurrence, those containing
exactly one of them, and those containing none of the two.

The contribution of a subset~$U$ in the first class to the
first summand on the left is
$N^{c(\alpha)-c(\alpha|_U)+f(\alpha|_U)-1}\e^{m-|U|}$.
Associate to the subset~$U\subset \{1,2,\dots,m\}$
the subset $U'\subset \{1,2,\dots,m-1\}$, which is obtained
from~$U$ by identifying $l$ and $l+1$ and diminishing
each permuted element to the right of~$l+1$ by~$1$.
We are going to show that the contribution of~$U$
to the first summand on the right coincides with
the contribution of~$U'$ to the first summand on the left.
If the elements~$l$ and~$l+1$ belong to the same cycle of~$\alpha|_U$
(or, which is the same, to the same cycle of~$\alpha$), then the number of cycles
in the first permutation on the right, as well as the
number of faces of the permutation
is the same as  the number of cycles
in the first permutation on the left, as well as the
number of faces of the permutation.
If the elements~$l$ and~$l+1$ belong to different cycles of~$\alpha|_U$,
then the number of cycles is one less than
$c(\alpha|_U)$, while the numbers of faces are preserved.
The argument for the
second summand is the same.

The second class of subsets splits into pairs:
a pair is obtained by adding either~$l$ or~$l+1$
to a subset of $\{1,2,\dots,,m\}\setminus \{l,l+1\}$.
The contribution of a subset~$U$ in the second class to the first summand on the left coincides with the contribution
of the subset  $U\Delta\{l,l+1\}$ to
the second summand  on
the left,  whence  the total contribution of the
pair  $U$ and $U\Delta\{l,l+1\}$ to the difference
on the left is~$0$.

The contribution of a subset~$U$ in the third class to the difference on the left also is~$0$.

The same  is  true for the contribution of a subset
$U'\subset \{1,2,\dots,m-1\}$ not containing~$l$ to the
difference in the right hand side.

The theorem is proved.

\subsection{Proof of Theorem~\ref{t-rsc}}

Pick an arc diagram~$B$ with~$n$  arcs
and represent it as an involution~$\alpha$ without fixed points
of the set $\{1,\dots,2n\}$.
Let~$U$ be a $k$-element subset of the set $V(B)$ of arcs in~$B$. Under the substitution  $u=2\e+N\e^2$, $w=N$
of the theorem, the contribution of the subset~$U$
to the refined skew-characteristic polynomial
$\overline{Q}_{\gamma(B)}$ becomes
$(2\e+N\e^2)^{n-k}N^{f(B|_U)-1}$.
We are going to show that this contribution coincides
with the total contribution to $F_\e(B)$ of the subsets
$U'\subset\{1,2,\dots,2n\}$ all whose elements are
ends of the arcs in~$V(B)\setminus U$, and for each arc in~$V(B)\setminus U$
at least one of its ends belongs to~$U'$. Note that any
subset~$U'$ of permuted points determines the subset~$U$
of arcs uniquely.

Pick a one-to one correspondence between the arcs
in~$V(B)\setminus U$
and the linear factors in $(2\e+N\e^2)^{n-k}$.
Inside each such factor, associate the term $N\e^2$ to
the subset of the ends consisting of both ends of the arc,
and pick a one-to-one correspondence between the two
ends of the arc on one side, and the two entrances of~$\e$
in the summand on the other. Then each subset
$U'\subset\{1,2,\dots,2n\}$ whose set of arcs is~$U$
determines uniquely a summand in each factor $2\e+N\e^2$.
The product of these summands times $N^{f(B|_{V(G)\setminus U})-1}$ coincides with the contribution of~$U'$ to $F_\e$. Indeed, each pair of
arc ends contributes $\e^2N$, since its presence in~$U'$
decreases the number of cycles by~$1$, while each
single end contributes the factor~$\e$, because
it does not affect the number of cycles in the permutation.

The theorem is proved.

\subsection{Proof of Theorem~\ref{t-is}}

The substitution $N=z^2-1$, $\e=\frac1{1-z}$ in the function~$F_\e$ leads to the
substitution $u=1,w=z^2-1$ in the refined skew characteristic polynomial $\overline{Q}_G(u,w)$.
Since the interlace polynomial $L_G(x)$ coincides with $\overline{Q}(1,x)$,
we obtain the assertion of the theorem.

\section{Odds and ends}

In this section, we present questions and results
relating the skew characteristic and the interlace polynomial
to Lie algebra weight systems and the Hopf algebra structures
on spaces of certain combinatorial objects.

\subsection{Hopf algebra structure on the spaces of graphs, chord diagrams,
and permutations}

Graphs, as well as chord diagrams and permutations, span a connected graded
commutative Hopf algebra. We refer the reader to~\cite{KL22}, for example.
Denote by~$\cG_n$ the vector space spanned by graphs with~$n$ vertices (over~$\C$, for definiteness).
Then the infinite direct sum
$$
\cG=\cG_0\oplus\cG_1\oplus\cG_2\oplus\dots
$$
is naturally endowed with the following two operations:
\begin{itemize}
    \item \emph{multiplication} $m:\cG\otimes\cG\to\cG$ induced by the disjoint union of graphs;
    \item \emph{comultiplication} $\mu:\cG\to\cG\otimes\cG$ induced by the operation splitting the set~$V(G)$ of vertices of a graph~$G$ into two disjoint subsets in all possible ways.
\end{itemize}
These two operations (together with natural unit, counit and antipode) make~$\cG$ into
a connected graded commutative cocommutative Hopf algebra.
These operations preserve $4$-term relations for graphs, which allows one to define
the corresponding quotient Hopf algebra.

A Hopf algebra structure on the space spanned by delta-matroids and their various versions
is defined in a similar way~\cite{LZ17}. It also admits a factorization modulo
$4$-term relations for graphs. A Hopf algebra structure on the space of chord diagrams
appears naturally only after $4$-term relations are factored out~\cite{K93}.

An element~$p$ of a Hopf algebra with a comultiplication~$\mu$
is said to be \emph{primitive} if its coproduct is
$\mu(p)=1\otimes p+p\otimes1$. According to the Milnor--Moore theorem, any connected
graded commutative cocommutative Hopf algebra is isomorphic to a polynomial Hopf algebra
in its primitive elements. In  particular, this is true for the Hopf algebra~$\cG$.
This implies that each homogeneous vector subspace~$\cG_n$ admits a decomposition into
a direct sum $\cG_n=P(\cG_n)\oplus D(\cG_n)$, Here $P(\cG_n)\subset\cG_n$ is the vector
subspace of primitive elements of degree~$n$, while $D(\cG_n)\subset\cG_n$ is the vector
subspace of decomposable elements, that is polynomials in primitive elements of smaller
degree (spanned, in turn, by disconnected graphs). This decomposition determines a projection~$\pi$ from the space~$\cG$ to the
space $P(\cG)$ of its primitive elements, which takes each primitive element to itself,
and each decomposable element to~$0$.

The space of linear functions on a Hopf algebra (or, more generally, on a coalgebra)
is endowed with a natural multiplication
called the convolution product. The \emph{convolution product} $f\circ g$ of two linear
functions~$f,g$ is the linear function defined as $f\circ g(a)=(f\otimes g)(\mu(a))$,
for any~$a\in\cG$. The projection~$\pi$ to the subspace of primitive elements
is nothing but the logarithm of the identity mapping ${\rm Id}:\cG\to\cG$,
where the logarithm is understood in terms of the convolution product,
see details in~\cite{KL22}, for example.

\begin{theorem}
For an arbitrary chord diagram~$B$ with at least two chords,
the value $F_\e(\pi(B))$ is independent of~$\e$.
\end{theorem}

The only exception is the chord diagram $K_1$ with a single chord, for which
$F_\e(K_1)=N+2\e+N\e^2=F_\e(\pi(K_1))$.

\begin{corollary}
For an arbitrary chord diagram~$c$ with at least two chords,
the value $F_\e(\pi(c))$ coincides with $F(\pi(c))$.
\end{corollary}

In other words, the family~$F_\e$ on chord diagrams can be defined as a
multiplicative weight system with values in the ring of
of polynomials in~$\e$ whose coefficients are polynomials in~$N$
such that their value on primitive elements coincides with that of~$F$,
with the exception of the chord diagram with a single chord, where it is $F_\e(K_1)$.

\subsection{Skew characteristic polynomial and $\gl(1|1)$-weight system}

The skew characteristic polynomial is the convolution
product of the following two
multiplicative graph invariants.
The first one takes a graph~$G$ to $N^{|V(G)|}$.
The second invariant takes the graph on a single
vertex to~$u$ and any other connected graph to~$0$.

The construction of Lie algebra weight systems can be extended to Lie superalgebras,
in particular, to the series $\gl(M|N)$.
It is shown in~\cite{ZY23} that for arbitrary~$M,N$
the $\gl(M|N)$-weight system can be induced from the universal $\gl$-weight
system. The $\gl(1|1)$-weight system is a special case of the $\gl(M|N)$-weight
system, for $M=N=1$.
For the $\gl(1|1)$-weight system, the corresponding substitution
looks like $C_0=M-N=0$, and the higher Casimir elements  $C_3,C_4,\dots$
can be expressed in terms of the first and the second Casimir
elements $C_1,C_2=$ by means of the generating function
$$
\sum_{k=1}^\infty C_kt^k
=\frac{C_1t}{\left(1-\frac{(-C_1^2+C_1+C_2)t}{2C_1}\right)\left(1-\frac{(C_1^2-C_1+C_2)t}{2C_1}\right)}.
$$
The following relation between the $\gl(1|1)$ weight system and the skew characteristic
polynomial has been remarked in~\cite{DL22}:

\begin{theorem}
    The skew characteristic polynomial
coincides with the result of substitution $C_1=1,C_2=u$ in the
$\gl(1|1)$-weight system expressed in terms
of the Casimir elements~$C_1,C_2$.
\end{theorem}

Note that under the substitution in the theorem the generating function for the
$\gl(1|1)$ Casimir elements acquires the form
\begin{eqnarray*}
\sum_{k=1}^\infty C_kt^k
&=&\frac{t}{\left(1-\frac{u}{2}t\right)^2}\\
&=&t+ut^2+\frac34u^2t^2+\frac12u^3t^3+\dots.
\end{eqnarray*}

It is easy to see
that the skew characteristic polynomial is the convolution product of the
    \emph{nondegeneracy} of a graph (the function
    taking on a graph value~$1$ if its adjacency matrix
    is nondegenerate over the two-element field,
    and value~$0$ otherwise)
    and the graph invariant taking
    on a graph~$G$ the value $u^{|V(G)|}$~\cite{DL22}. The refined skew characteristic polynomial
    is the convolution product of the
    skew characteristic polynomial and the graph invariant taking
    on a graph~$G$ the value $N^{{\rm corank}~A_G}$.

It follows immediately that the value of the skew characteristic polynomial on the projection~$\pi(G)$    of a graph~$G$ to the subspace of primitives is a constant. Similarly,
     the value of the refined skew characteristic polynomial on the projection~$\pi(G)$
    of a graph~$G$ to the subspace of primitives is the same constant on all graphs
    with at least two vertices; on the graph~$K_1$ having a single vertex
    this value is $1+N$.

Indeed, the projection to the subspace of primitives is the logarithm (with respect to the convolution
product) of the identity mapping of the Hopf algebra of graphs to itself.
Applying the logarithm to the convolution product of two invariants, we obtain the sum of the logarithms of the factors.
The invariant taking a graph~$G$ to~$u^{|V(G)|}$, when projected to the subspace of primitives,
is~$u$ for the graph~$K_1$ and~$0$ for any graph with more than one vertex.

\begin{lemma}
    For intersection graphs of chord diagrams,
    the  graph invariant equal to $N^{|V(G)|}$  on  a graph~$G$
    is induced from the $\gl$-weight system under the substitution $C_m=N$,
    $m=1,2,3,\dots$.
\end{lemma}

Indeed, the substitution in the lemma takes each permutation~$\alpha$ to~$N^{c(\alpha)}$,
and the number of cycles in a permutation defined by an arc diagram coincides
with the number of arcs in it, that is, with the number of vertices in its intersection
graph.

\subsection{Interlace polynomial for permutations and partial duality}

The  substitution in the  family $F_\e$, which makes it,
for arc diagrams, into the interlace polynomial,
can be done for  arbitrary permutations, not only for
involutions without fixed points. In spite  of the fact
that after this substitution  we obtain,  generally
speaking, a rational function rather than a polynomial,
we define the \emph{interlace polynomial $L(\alpha)$ of a permutation $\alpha\in\S_m$}
as the function
$$
L(\alpha)=F_\e(\alpha)|_{N=z^2-1,\e=\frac1{1-z}}.
$$
This rational  function either has a unique pole at the point $z=1$, or is a polynomial.

For example, $F_\e((1,2,3))=N\e^3+3\e^2+3N\e +N^2$, whence
$$
L((1,2,3))=\frac{z^2(z-2)(z^3-z-1)}{(z-1)^2},
$$
while
$F_\e((1,3,2))=N\e^3+3\e^2+3N\e +1$, whence
$$
L((1,3,2))=-\frac{z^2(3z-4)}{(z-1)^2}.
$$

It happens that in all the examples we  have
computed, the power series expansions of the
rational functions $L(\alpha)$ at $z=0$
have positive coefficients. In particular,
in the examples above
$$
\frac{z^2(z-2)(z^3-z-1)}{(z-1)^2}=
2z^2+5z^3+7z^4+7z^5+8z^6+9z^7+\dots
$$
and
$$
-\frac{z^2(3z-4)}{(z-1)^2}=
4z^2+5z^3+6z^4+7z^5+8z^6+9z^7+\dots
$$
It would be interesting to know whether this
assertion is true for arbitrary permutations,
and if yes, then what is the
combinatorial meaning  of the coefficients of the corresponding power series.

Now let $\alpha\in\S_m$ be a permutation, and suppose
$a,b\subset\{1,2,\dots,m\}$ are two intersecting
orbits of~$\alpha$ of length~$2$ (that is, the elements
in these two orbits alternate). Then we can define the
\emph{pivot} $\alpha^{ab}\in\S_m$ of~$\alpha$ in the same way as
in the case of arc diagrams (which is exactly the partial duality for hypermaps~\cite{CV22}). We will  give a proof of the following
two theorems elsewhere.

\begin{theorem}
    The interlace polynomial of permutations is invariant under pivot,
$$
L(\alpha)=L(\alpha^{ab})
$$
for any permutation~$\alpha$ and any two its intersecting
orbits~$a,b$.
\end{theorem}

For example, the permutation $(6,5,8,7,2,3,4,9,1)$ is the pivot of the permutation $(3,5,6,7,2,8,4,9,1)$ with respect to the pair of
the chords $a=(2,5)$ and $b=(4,7)$. The value of the
interlace polynomial on each of the two permutations is
$$
\frac{z^2 (6 - 2 z + 7 z^2 + z^3 + z^4 - 8 z^5 - 3 z^6 - z^7 + 10 z^8 -
   6 z^9 + z^{10})}{(1 - z)^4}.
$$

For permutations having two intersecting chords (orbits of length~$2$),
we can express the value of the interlace
polynomial in simpler permutation in a way
similar to the recurrence relations~(\ref{eIPr}) in the definition
of the interlace polynomial for graphs.

\begin{theorem}
    For a permutation $\alpha$, the following relation holds:
    $$L(\alpha)= L(\alpha|_{\{1,\dots,m\}\setminus\{a,b\}}) + L(\alpha^{ab}|_{\{1,\dots,m\}\setminus\{a,b\}}).
    $$
   \end{theorem}

Behavior of the interlace polynomial with respect to generalized duality for
hypermaps~\cite{CV22} also deserves being studied.

\end{document}